\documentclass{amsart}
\usepackage[dvips,final]{graphics}
\usepackage{array}
\usepackage{arydshln}
\usepackage[makeroom]{cancel}
 \usepackage[all]{xy}
 \usepackage{url}
\usepackage{multirow, blkarray}
\usepackage{booktabs}
\usepackage{textcomp}
 \usepackage[final]{epsfig}
 \usepackage{color}
\usepackage[T1]{fontenc}      
\usepackage[english,french]{babel}
\usepackage[utf8]{inputenc}
\usepackage{blindtext}

\usepackage{amsfonts,amscd,array, mathdots, epigraph}
\usepackage{amsmath}
\usepackage{amssymb}
\usepackage{amsthm}
\usepackage{mathrsfs}
\usepackage{stmaryrd}
\usepackage{slashbox}
\usepackage{diagbox}
\usepackage{enumitem}

\usepackage{ulem}
\usepackage{tikz}
\usepackage{xcolor}
\usepackage{multicol}
\usepackage{tcolorbox}
\definecolor{ufogreen}{rgb}{0.24, 0.82, 0.44}

\usepackage{hyperref}

\begin{document}

% THEOREMS -------------------------------------------------------

\newtheorem{theorem}{Theorem}[section]
\newtheorem{theore}{Theorem}
\newtheorem{definition}[theorem]{Definition}
\newtheorem{proposition}[theorem]{Proposition}
\newtheorem{corollary}[theorem]{Corollary}
\newtheorem{con}{Conjecture}
\newtheorem*{remark}{Remark}
\newtheorem*{remarks}{Remarks}
\newtheorem*{pro}{Probleme}
\newtheorem*{examples}{Examples}
\newtheorem*{example}{Example}
\newtheorem{lemma}[theorem]{Lemma}

\newtcolorbox{attentionbox}{
  colback=yellow!10,
  colframe=red!70!black,
  title=Warning,
  fonttitle=\bfseries
}

%---------------------------------------------

\title{Maximal size of irreducible $\lambda$-quiddities over polynomial and formal power series rings}

\author{Flavien Mabilat}

\date{}

\email{flavien.mabilat@univ-reims.fr}

\maketitle

\selectlanguage{french}

\begin{abstract}

L'étude de la combinatoire du groupe modulaire et des frises de Coxeter amènent naturellement à étudier une équation matricielle, parfois nommée équation de Conway-Coxeter, dont les solutions de taille $n$, appelés $\lambda$-quiddités, sont des $n$-uplets d'éléments de l'anneau $B$ considéré. La bonne connaissance de ces dernières repose sur une notion de solutions irréductibles qui permettent de reconstruire l'ensemble des $n$-uplets recherchés. L'un des problèmes centraux qui émerge alors naturellement de cette démarche est de savoir si les $\lambda$-quiddités irréductibles sur $B$ ont une taille bornée et d'obtenir le cas échéant une telle borne. Ici, on se propose de lister des résultats qui répondent à cette question dans le cas des anneaux de polynômes $A[X]$ et $\mathbb{K}[X]$, avec $A$ un anneau commutatif unitaire fini et $\mathbb{K}$ un corps commutatif. Par ailleurs, les résultats énoncés permettront également de considérer facilement de nombreuses situations dans lesquelles $A$ est infini. Pour finir, on donnera la réponse complète à la question initiale pour tous les anneaux de séries formelles.
\\
\end{abstract}

\selectlanguage{english}
\begin{abstract}

The study of the combinatorics of the modular group and of Coxeter's friezes naturally leads to the investigation of a matrix equation, sometimes referred to as the Conway–Coxeter equation. The solutions of size $n$ of this equation, called $\lambda$-quiddities, are $n$-tuples of elements of a given ring $B$. A detailled understanding of these objects relies on the notion of irreducible solutions, from which all $\lambda$-quiddities can be reconstructed. One of the central questions that naturally arises in this context is whether the irreducible $\lambda$-quiddities over $B$ have bounded size, and, if so, how to determine such a bound. In this paper, we aim to list results that address this question in the case of polynomial rings $A[X]$ and $\mathbb{K}[X]$, where $A$ is a finite commutative unitary ring and $\mathbb{K}$ is a commutative field. Moreover, the stated results will also make it possible to treat easily many situations in which $A$ is infinite. Finally, we shall give a complete answer to the initial question for all rings of formal power series.
\\
\end{abstract}

\selectlanguage{french}

\thispagestyle{empty}

\begin{attentionbox}
This document does not contain any proofs. Its purpose is to make the results presented in the French article \cite{Mprin}, as well as the motivations that led to obtaining them, easily accessible to users of arXiv.
\end{attentionbox}

\section{Introduction}
\label{intro}

For several years, the study of Coxeter's friezes—initially introduced in the early 1970s by H. S. M. Coxeter—has experienced a strong resurgence of interest, and numerous papers have been devoted to these objects, which appear simple at first glance but possess fascinating properties. Originally conceived as tools for the study of the Pentagramma mirificum (see \cite{Cox}), they have since become a subject of study in their own right, extremely active, with countless ramifications and generalizations (see \cite{Mo1}). More precisely, friezes are arrangements of numbers in the plane belonging to a fixed set $R$ and satisfying an arithmetic rule, called the unimodular rule, to which one often adds a property known as tameness (see \cite{CH,Mo1} for precise definitions). When $R$ is a commutative unitary ring, the construction of tame Coxeter's friezes is directly related to the solutions of the following matrix equation (see \cite{CH}, Proposition 2.4):
\[M_{n}(a_1,\ldots,a_n):=\begin{pmatrix}
    a_{n} & -1_{R} \\
    1_{R}  & 0_{R} 
   \end{pmatrix} \ldots \begin{pmatrix}
    a_{1} & -1_{R} \\
    1_{R}  & 0_{R} 
   \end{pmatrix}=-Id .\]

The occurrence of the matrices $M_{n}(a_1,\ldots,a_n)$ is all the more interesting since they already arise in a myriad of seemingly very different contexts. In particular, they appear in the study of convergents of Hirzebruch-Jung continued fractions, as well as in the matrix formalization of discrete Sturm-Liouville equations. However, their most notable use is in the representation of elements of the modular group $SL_{2}(\mathbb{Z})$. Indeed, by using the standard generators of this group,
\[T=\begin{pmatrix}
 1 & 1 \\[2pt]
    0    & 1 
   \end{pmatrix}~{\rm and}~S=\begin{pmatrix}
   0 & -1 \\[2pt]
    1    & 0 
   \end{pmatrix},
 \] 
\noindent one can show that, for every matrix $M$ in the modular group, there exist strictly positive integers $n$ and $a_{1},\ldots,a_{n}$ such that
$M=M_{n}(a_1,\ldots,a_n)=T^{a_{n}}S\ldots T^{a_{1}}S$. Since this factorization is not unique, one is naturally led to seek all possible factorizations of a given matrix or set of matrices. In particular, a complete solution to this problem is available in the cases where $M=\pm Id$ (see \cite{O}), $M=\pm S$ and $M=\pm T$ (see \cite{M9}), as well as numerous partial results for congruence subgroups (see in particular \cite{M1,M4}). 
\\
\\\indent All of this has led several authors, notably V. Ovsienko and M. Cuntz, to investigate a natural generalization of the equation used in the study of Coxeter's friezes (see \cite{C,O}). Specifically, one fixes a set $A$—which, in this text, will be a commutative unitary ring—and considers the following matrix equation:
\begin{equation}
\label{p}
\tag{$E_{A}$}
M_{n}(a_1,\ldots,a_n)=\pm Id.
\end{equation}

\noindent This equation is sometimes called the Conway–Coxeter equation, and its $n$-tuple solutions are referred to as $\lambda$-quiddities of size $n$ over $A$. Note that the search for various factorizations of elements in the congruence subgroups mentioned in the previous paragraph is equivalent to solving \eqref{p} over the rings $\mathbb{Z}/N\mathbb{Z}$.
\\
\\\indent There are several approaches to the study of $\lambda$-quiddities over a given set. One may seek a geometric construction of all the solutions (see \cite{O}), attempt to determine rings for which there are only finitely many $\lambda$-quiddities of fixed size (see \cite{CHP}), or count the solutions of fixed size over a finite set (see, for example, \cite{BC,CO2,CM,Mo2}). However, the potential of these approaches remains limited, as there are few situations in which substantial results can be obtained. In most cases, the most effective tool is a notion of irreducibility introduced by M. Cuntz (see \cite{C} and the following section). This notion indeed allows one to restrict the study of \eqref{p} to a thorough understanding of the irreducible $\lambda$-quiddities over $A$, since the latter can be used to reconstruct all solutions. The most convenient cases are, of course, those in which a complete list of irreducible $\lambda$-quiddities can be obtained. In the majority of cases, however, this task proves too complex due to the sheer, sometimes overwhelming, number of irreducible solutions and the diversity of their forms. Moreover, compiling excessively long lists of irreducibles is of extremely limited value, even if complete. Thus, one usually restricts oneself to exhaustively listing the irreducibles only for rings of small cardinality (see \cite{M1}, Section 4, and \cite{M7}, Theorems 2.6 and 2.7). For other rings, the strategy is rather to construct simple families of $\lambda$-quiddities with good irreducibility properties and to investigate general properties satisfied by all irreducible solutions. The most important of these properties is the possible existence of a finite upper bound on the size of irreducibles.
\\
\\\indent On this topic, the most important result is the existence of such a bound for all finite rings (see \cite{M} and Section \ref{ldea}). However, when considering an infinite ring $A$, satisfactory results are available only in a very limited number of cases (see Theorem \ref{32}). To partially address this gap, we state here, without proofs, results concerning three highly important families of infinite rings: polynomial rings, formal power series rings, and infinite fields. The proofs of these results can be found in \cite{Mprin}.

\section{Definitions}
\label{prin}

The aim of this section is to give the definition of the objects mentionned in the introduction and used in the next section.
\\
\\\indent In this text, all rings are commutative and unitary and different from $\{0_{A}\}$. Let $(A,+,\times)$ be such a ring. ${\rm car}(A)$ is the characteristic of $A$. $0_{A}$ denotes the additive identity and $1_{A}$ the multiplicative identity. We denote by $U(A)$ the group of units of $A$, and by $A[[X]]$ the ring of formal power series over $A$. If $A$ and $B$ are isomorphic rings, we write $A \cong B$. We denote by $\mathbb{N}$ the set of non-negative integers and $\mathbb{N}^{*}$ the set of positive integers. If $N \geq 2$ we set, for all integers $a$, $\overline{a}:=a+N\mathbb{Z}$.
\\
\\We begin by the precise definition of $\lambda$-quiddities.

\begin{definition}[\cite{C}, Definition 2.2]
\label{21}

Let $n \in \mathbb{N}^{*}$. The $n$-tuple $(a_{1},\ldots,a_{n}) \in A^{n}$ is a $\lambda$-quiddity over $A$ if $(a_{1},\ldots,a_{n})$ is a solution of \eqref{p}, that is to say if there exists $\epsilon \in \{\pm 1_{A}\}$ such that $M_{n}(a_{1},\ldots,a_{n})=\epsilon Id$. If there is no ambiguity, we will omit to precise the ring. 

\end{definition}

Throughout the rest of this text, we will speak indiscriminately of the solutions of \eqref{p} or of the $\lambda$-quiddities over $A$. To define the notion of irreducibility, we need the two following definitions.

\begin{definition}[\cite{C}, Lemma 2.7]
\label{22}

Let $(a_{1},\ldots,a_{n}) \in A^{n}$ and $(b_{1},\ldots,b_{m}) \in A^{m}$. We define the following operation : \[(a_{1},\ldots,a_{n}) \oplus (b_{1},\ldots,b_{m})= (a_{1}+b_{m},a_{2},\ldots,a_{n-1},a_{n}+b_{1},b_{2},\ldots,b_{m-1}).\] The $(n+m-2)$-tuple obtained is the sum of $(a_{1},\ldots,a_{n})$ with $(b_{1},\ldots,b_{m})$.

\end{definition}

The operation $\oplus$ is very useful since it has the interesting following property. Let $(b_{1},\ldots,b_{m})$ be a solution of \eqref{p}. $(a_{1},\ldots,a_{n}) \oplus (b_{1},\ldots,b_{m})$ is a $\lambda$-quiddity over $A$ if and only if $(a_{1},\ldots,a_{n})$ is a solution of \eqref{p} (see \cite{C} Lemma 2.7 and \cite{WZ} Lemma 1.9).

\begin{definition}[\cite{C}, Definition 2.5]
\label{23}

Let $(a_{1},\ldots,a_{n}) \in A^{n}$ and $(b_{1},\ldots,b_{n}) \in A^{n}$. The $n$-tuple $(a_{1},\ldots,a_{n})$ is said to be equivalent to $(b_{1},\ldots,b_{n})$ (denoted by $(a_{1},\ldots,a_{n}) \sim (b_{1},\ldots,b_{n})$) if $(b_{1},\ldots,b_{n})$ can be obtained by cyclic rotations of $(a_{1},\ldots,a_{n})$ or of $(a_{n},\ldots,a_{1})$.

\end{definition}

The relation $\sim$ is an equivalence relation on the set $A^{n}$ (see \cite{WZ} Lemma 1.7). Moreover, if $(a_{1},\ldots,a_{n}) \sim (b_{1},\ldots,b_{n})$ then $(a_{1},\ldots,a_{n})$ is a solution of \eqref{p} if and only if $(b_{1},\ldots,b_{n})$ is a $\lambda$-quiddity over $A$ (see \cite{C} Proposition 2.6).
\\
\\We can now define the notion of irreducibility mentioned in the introduction.

\begin{definition}[\cite{C}, Definition 2.9]
\label{24}

A $\lambda$-quiddity over $A$ $(c_{1},\ldots,c_{n})$ ($n \geq 3$) is said to be reducible if there exists a $\lambda$-quiddity over $A$ $(b_{1},\ldots,b_{l})$ and an $m$-tuple $(a_{1},\ldots,a_{m}) \in A^{m}$ such that \begin{itemize}
\item $(c_{1},\ldots,c_{n}) \sim (a_{1},\ldots,a_{m}) \oplus (b_{1},\ldots,b_{l})$;
\item $m \geq 3$ and $l \geq 3$.
\end{itemize}
A $\lambda$-quiddity over $A$ is said to be irreducible if it is not reducible. We denote by $Irr(A)$ the set of irreducible $\lambda$-quiddities over $A$.

\end{definition}

\begin{remark} 

{\rm The 2-tuple $(0_{A},0_{A})$ is never considered as an irreducible solution of \eqref{p}.}

\end{remark}

As announced in the introduction, this notion of irreducibility allows us to restrict the study of \eqref{p} to a thorough understanding of the irreducible $\lambda$-quiddities over $A$, since the latter can be used to reconstruct all solutions. Since compiling a complete list of the elements of ${\rm Irr}(A)$ for a fixed ring is generally out of reach, we focus on identifying properties satisfied by irreducible $\lambda$-quiddities. To this end, we make particular use of the object introduced in the following definition:

\begin{definition}
\label{25}

Let $A$ be a commutative and unitary ring. If the size of the irreducible solutions of \eqref{p} are bounded, we set $\ell_{A}$ to be the maximal size of them. In the other case, we set $\ell_{A}=+\infty$.

\end{definition}

The aim of this paper is to give, without proofs, the exact values of $\ell_{A}$ for the infinite commutative fields, the polynomial rings over finite rings and the rings of formal power series.

\section{Determination of $\ell_{A}$ for certain rings}

\noindent The aim is to section is to give all the known values of $\ell_{A}$.

\subsection{Some already kwown results}
\label{ldea}

When studying irreducible $\lambda$-quiddités over a ring $A$, two questions naturally arise. Is there a finite number of irreducible solutions over $A$ ? Are their sizes bounded ? In the case of finite rings, these two questions are in fact equivalent, and the common answer to both is provided by the result below.

\begin{theorem}[\cite{M}, Theorem 1.1]
\label{31}

Let $(A,+,\times)$ be a finite commutative and unitary ring.
\begin{itemize}
\item If $car(A)=2$ then $4 \leq \ell_{A} \leq \frac{\left|SL_{2}(A)\right|}{\left|A\right|}+2$.
\item If $car(A) \neq 2$ then ${\rm max}(4,car(A)) \leq \ell_{A} \leq \frac{\left|SL_{2}(A)\right|}{2\left|A\right|}+2$.
\end{itemize}

\end{theorem}

Note that the general lower bounds provided by this theorem are optimal. Indeed, $\ell_{\mathbb{Z}/2\mathbb{Z}}=4$ and $\ell_{\mathbb{Z}/6\mathbb{Z}}=6$. For rings of small cardinality, the exact values are known (see, for instance, \cite{CM}, Section 5.2).

\medskip

When the ring under consideration is infinite, each case must be studied separately. We collect in the theorem below several known values

\begin{theorem}
\label{32}

i) (\cite{CH}, Theorem 6.2) $\ell_{\mathbb{Z}}=4$. Moreover,
\[Irr(\mathbb{Z})=\{\{(1,1,1), (-1,-1,-1), (0,a,0,-a), (a,0,-a,0); a \in \mathbb{Z}-\{\pm 1\} \}.\]

\noindent ii) (\cite{C}, Proposition 3.3) $\ell_{\mathbb{Z}[i]}=+\infty$.
\\
\\iii) (\cite{M2}, Theorem 2.7) Let $\alpha$ be a transcendental complex number. We have $\ell_{\mathbb{Z}[\alpha]}=4$ and :
\[Irr(\mathbb{Z}[\alpha])=\{(1,1,1), (-1,-1,-1), (0,P(\alpha),0,-P(\alpha)), (P(\alpha),0,-P(\alpha),0); P \in \mathbb{Z}[X]-\{\pm 1\} \}.\]

\noindent iv) (\cite{M7}, Theorem 2.8) Let $(A_{i})_{i \in I}$ be a family of commutative rings. We assume that at least two rings in this family have characteristic 
0. We have $\ell_{\prod_{i \in I} A_{i}}=+\infty$. In particular, $\ell_{\mathbb{Z} \times \mathbb{Z}}=+\infty$.

\end{theorem}

\begin{proposition}[\cite{M7}, Corollary 3.3]
\label{343}

Let $A$ and $B$ be isomorphic commutative unitary rings, and let $f:A\rightarrow B$ be an isomorphism. An $n$-tuple $(a_{1},\ldots,a_{n})$ is an irreducible $\lambda$-quiddity over $A$ if and only if $(f(a_{1}),\ldots,f(a_{n}))$ is an irreducible $\lambda$-quiddity over $B$. In particular, $\ell_{A}=\ell_{B}$.

\end{proposition}

\begin{proposition}[\cite{M}, Proposition 2.16]
\label{342}

Let $A$ and $B$ be two unitary commutative rings. We suppose $A \subset B$ and $1_{A}=1_{B}$. Any irreducible solution over $A$ is an irreducible solution over $B$. In particular, one has $\ell_{A} \leq \ell_{B}$.

\end{proposition}

\noindent The following result is an immediate consequence of part (iii) of Theorem~\ref{32} and Proposition~\ref{343}.

\begin{corollary}[\cite{Mprin}, Corollary 3.7]
\label{341}

We have $\ell_{\mathbb{Z}[X]}=4$. Moreover, 
\[Irr(\mathbb{Z}[X])=\{\{(1,1,1), (-1,-1,-1), (0,P(X),0,-P(X)), (P(X),0,-P(X),0); P \in \mathbb{Z}[X]-\{\pm 1\} \}.\]

\end{corollary}

\subsection{Values of $\ell_{A}$ over polynomial rings}

We give here all the recent results concerning the concrete value of $\ell_{A[X]}$.

\begin{theorem}[\cite{Mprin}, Theorem 2.6]
\label{27}

Let $A$ be a commutative and unitary finite ring.
\[\ell_{A[X]}<+\infty \Longleftrightarrow A\cong \mathbb{Z}/2\mathbb{Z}~or~A \cong \mathbb{Z}/3\mathbb{Z}.\]
\noindent Moreover, $\ell_{(\mathbb{Z}/2\mathbb{Z})[X]}=\ell_{(\mathbb{Z}/3\mathbb{Z})[X]}=4$ and, for $k \in \{2,3\}$, we have :
\[Irr((\mathbb{Z}/k\mathbb{Z})[X])=\{\pm (\overline{1},\overline{1},\overline{1}), (\overline{0},P(X),\overline{0},-P(X)), (P(X),\overline{0},-P(X),\overline{0}); P \in (\mathbb{Z}/k\mathbb{Z})[X]-\{\pm \overline{1}\}\}.\]

\end{theorem}

\begin{proposition}[\cite{Mprin}, Proposition 2.7]
\label{26}

Let $A$ be a commutative unitary ring and $n \geq 2$.
\\
\\i) If $A$ is finite then : $\ell_{A[X_{1},\ldots,X_{n}]}<+\infty \Longleftrightarrow A \cong \mathbb{Z}/2\mathbb{Z}~or~A \cong \mathbb{Z}/3\mathbb{Z}$.
\\
\\ii) We assume that $A$ satisfies at least one of the four properties below:
\begin{itemize}
\item $car(A) \notin \{0,2,3\}$;
\item $A$ contains a nilpotent element;
\item $A$ is decomposable;
\item $A$ contains an invertible element different from $\pm 1_{A}$. 
\end{itemize}

\noindent We have $\ell_{A[X_{1}]}=\ell_{A[X_{1},\ldots,X_{n}]}=+\infty$.

\end{proposition}

\begin{examples}
{\rm 
\begin{itemize}
\item $A:=\frac{\mathbb{Z}[Y]}{<Y^{n}>}$, with $n \geq 2$. Since $(Y+<Y^{n}>)^{n}=0_{A}$, $\ell_{A[X]}=+\infty$.
\item $A:=\frac{\mathbb{Z}[Y]}{<Y^{2}+Y+1>}$. Since we have $(-Y+<Y^{2}+Y+1>)(Y+1+<Y^{2}+Y+1>)=1_{A}$ and $(-Y+<Y^{2}+Y+1>) \neq \pm 1_{A}$, $\ell_{A[X]}=+\infty$.
\item $A:=\mathcal{F}(\mathbb{Z},\mathbb{Z}/4\mathbb{Z})$. Since ${\rm car}(A)=4$, $\ell_{A[X]}=+\infty$.
\item $A:=\mathbb{Z} \times (\mathbb{Z}/2\mathbb{Z})$. Since $A$ is the direct product of a finite number of rings, $\ell_{A[X]}=+\infty$.
\\
\end{itemize}
}
\end{examples}

\begin{proposition}[\cite{Mprin}, Proposition 4.9]

Let $n \geq 1$. $\ell_{\mathbb{Z}[X_{1},\ldots,X_{n}]}=4$.

\end{proposition}

\begin{proposition}[\cite{Mprin}, Proposition 4.10]

Let $k \geq 4$ and $n \geq 1$. $\ell_{\mathbb{Z}[i\sqrt{k}]}=\ell_{(\mathbb{Z}[i\sqrt{k}])[X_{1},\ldots,X_{n}]}=4$.

\end{proposition}

In this section, we have seen that if $A$ is an infinite ring verifying $\ell_{A[X]}<+\infty$ then $A$ is a reduced, indecomposable ring with at most two units. One might be led to believe that this condition is sufficient. In fact, this is not the case. Consider, for instance,  $A:=\mathbb{Z}[i\sqrt{2}]$. The ring $A$ is an integral domain, hence in particular indecomposable and reduced, and it has only two units, namely $1$ and $-1$. Inspired by the solution constructed for $\mathbb{Z}[i]$ (see \cite{C}, Proposition 3.3), one may consider, for every $l \geq 1$, the following tuple:

\[(i\sqrt{2},\sqrt{2}-i\sqrt{2},\underbrace{\sqrt{2},2\sqrt{2},\ldots,\sqrt{2},2\sqrt{2},\sqrt{2}}_{2l+1},\sqrt{2}+i\sqrt{2},-i\sqrt{2},\]
\[-\sqrt{2}+i\sqrt{2},\underbrace{-\sqrt{2},-2\sqrt{2},\ldots,-\sqrt{2},-2\sqrt{2},-\sqrt{2}}_{2l+1},-\sqrt{2}-i\sqrt{2}).\]

\noindent One shows that this $(4l+8)$-tuple is an irreducible solution over $A$. Consequently, $\ell_{A}=+\infty$, and therefore $\ell_{A[X]}=+\infty$.

\subsection{Values of $\ell_{A}$ over commutative fields}

We give here all the recent results concerning the concrete value of $\ell_{\mathbb{K}}$, with $\mathbb{K}$ a commutative field.
\\
\\\indent The results of the previous section readily yield that $\ell_{\mathbb{K}[X]}=+\infty$ when $\mathbb{K}$ is an infinite field. Indeed, an infinite field contains an invertible element distinct from $\pm 1_{\mathbb{K}}$, which allows one to apply Proposition \ref{26}. However, knowledge of $\ell_{\mathbb{K}[X]}$ only becomes truly meaningful once $\ell_{\mathbb{K}}$ is understood. Moreover, fields are rings of particular interest in their own right. It is therefore natural to seek a better understanding of the value of $\ell_{\mathbb{K}}$. This is the problem addressed in this subsection part. More precisely, we will give the general statement \og
$\ell_{\mathbb{K}}=+\infty$ when $\mathbb{K}$ is an infinite field.\fg. However, in the case of fields of characteristic 0, there exist irreducible $\lambda$-quiddités of arbitrarily large size that are particularly simple. We therefore begin by treating this special case separately.

\begin{proposition}[\cite{Mprin}, Proposition 5.2]

Let $n \geq 2$. The $(n+3)$-tuple $\left(2n+1,\frac{n+1}{2n+1},3,2,\ldots,2,\frac{2n}{2n+1}\right)$ is an irreducible $\lambda$-quiddity over $\mathbb{Q}$.

\end{proposition}

\begin{example}
{\rm Cyclotomic fields contain $\mathbb{Q}$. Consequently, by Proposition \ref{342}, the solution given in the previous proposition is irreducible over these fields.
}
\end{example}

\noindent We can draw inspiration from this construction to obtain the result below.

\begin{proposition}[\cite{Mprin}, Proposition 5.3]
\label{58}

Let $a \neq 0$ and $b \geq 2$ be two coprime integers, and let $\mathbb{D}$ denote the ring of decimal numbers. Set $A:=\mathbb{Z}\left[\frac{a}{b}\right]$. Then, $\ell_{A}=+\infty$. In particular, $\ell_{\mathbb{D}}=+\infty$.

\end{proposition}

\noindent In the general case, we have the following theorem.

\begin{theorem}[\cite{Mprin}, Theorem 2.8 and Proposition 5.5]
\label{54}

Let $\mathbb{K}$ be a commutative field. Let $n \geq 3$ be an odd integer. There exists an irreducible $\lambda$-quiddity of size $n$ over $\mathbb{K}$. In particular, $\ell_{\mathbb{K}}=+\infty$.

\end{theorem}

\begin{example}
{\rm We have $\ell_{(\mathbb{Z}/2\mathbb{Z})[t]}=4$ but $\ell_{(\mathbb{Z}/2\mathbb{Z})(t)}=+\infty$.
}
\end{example}

Before leaving the setting of fields, let us note that the results obtained for them make it possible to consider a whole range of rings that are neither fields nor polynomial rings. For instance, $\mathbb{Q}[\sqrt{2}]$ is a unitary commutative ring containing $\mathbb{Q}$. Hence, one has $\ell_{\mathbb{Q}[\sqrt{2}]}=+\infty$.

\subsection{The case of formal power series rings}

We give here the exact value of $\ell_{A[[X]]}$ for all commutative unitary rings $A$.

\begin{theorem}[\cite{Mprin}, Theorem 2.9]

Let $A$ be a commutative unitary ring. We have $\ell_{A[[X]]}=+\infty$.

\end{theorem}

\begin{example}
{\rm We have $\ell_{\mathbb{Z}[X]}=4$ but $\ell_{\mathbb{Z}[[X]]}=+\infty$.
}
\end{example}

\begin{corollary}[\cite{Mprin}, Corollary 5.7]

Let $A$ be a commutative unitary ring and $n \geq 1$. We have $\ell_{A[[X_{1},\ldots,X_{n}]]}=+\infty$.

\end{corollary}

\end{document}